\documentclass[leqno,10pt,twoside,bezier]{article}
\usepackage{amssymb,amsmath,amsfonts,amsthm,latexsym,ifthen,makeidx,bezier,epic,eepic}

       \newcommand{\kla}[1]{ {\langle #1 \rangle} }

       \newcommand{\st}{\;|\;}  
        \newcommand{\dom}{ {\rm dom}
       } \newcommand{\ran}{ {\rm ran} }

        \newcommand{\sub}{\subseteq}
        
        \newfont{\ssi}{cmssi12 at 12pt}

\newcommand{\rest}{\restriction} \newcommand{\cf}{ {\rm cf} }
 
\newcommand{\verl}{{{}^\frown}}
 
 \newcommand{\id}{ {\rm id} }
\newcommand{\qedd}[1]{\nopagebreak\hspace*{\fill}$ \Box_{#1} $\\}
\renewcommand{\qed}{\nopagebreak\hspace*{\fill}$ \Box $\\}

\newenvironment{ea*}{\begin{eqnarray*}}{\end{eqnarray*}}
\newtheorem{thm}{Theorem}[section]
\newtheorem{cor}[thm]{Corollary}
\newtheorem{lem}[thm]{Lemma}

\newtheorem{sealinglemma}[thm]{Sealing Lemma}

\newtheorem{question}[thm]{Question}
\newtheorem{statement}[thm]{Statement}
\newtheorem{mainthm}[thm]{Main Theorem}

\newtheorem*{thmStar}{Theorem}

\setbox0=\hbox{$\longrightarrow$} \newcommand{\To}{\longrightarrow}

\renewcommand{\proof}{\noindent{\em Proof. }}
\newcommand{\prooff}[1]{\vspace{12pt}\noindent{\em Proof of
    {#1}}.}

\newcommand{\vb}{{\vec{b}}} \newcommand{\vc}{{\vec{c}}}

\newcommand{\vs}{{\vec{s}}} \newcommand{\vt}{{\vec{t}}}
\newcommand{\vu}{{\vec{u}}} 
\newcommand{\vw}{{\vec{w}}} \newcommand{\vx}{{\vec{x}}}

\newcommand{\vC}{{\vec{C}}} 
\newcommand{\vE}{{\vec{E}}} 
 \newcommand{\vT}{{\vec{T}}}

 \renewcommand{\phi}{\varphi}
\newcommand{\card}{{\mbox{\rm card}}} \newcommand{\tx}[1]{\;\mbox{\rm
    #1}\;}  
 \newcommand{\ZFC}{{\tt ZFC}}
\newcommand{\ttt}{{\tt t}}

 \newcommand{\forces}{\Vdash}
 \newcommand{\Aut}{\mathop{\hbox{Aut}}}
\def\<#1>{\langle#1\rangle} \newcommand{\restrict}{\upharpoonright}
\renewcommand{\P}{\mathord{\Bbb P}}
\newcommand{\Q}{\mathord{\Bbb Q}}
\newcommand{\of}{\subseteq}
\newcommand{\intersect}{\cap} 
 \newcommand{\cross}{\times}
 \newcommand{\CF}{{\mbox{\rm CF}}}

\renewcommand{\diamond}{\diamondsuit}
\newcommand{\keinr}{\hspace*{-\parindent}}

\begin{document}
\hyphenation{covering}
\title{Changing the Heights of Automorphism Towers by Forcing with
  Souslin Trees over L}
\author{Gunter Fuchs\\Westf\"{a}lische Wilhelms-Universit\"{a}t M\"{u}nster  \and
  Joel David Hamkins\\
  The City University of New York}
\maketitle
\begin{abstract}
We prove that there are groups in the constructible universe whose automorphism towers are highly malleable by forcing. This is a consequence of the
fact that, under a suitable diamond hypothesis, there are sufficiently many highly rigid non-isomorphic Souslin trees whose isomorphism relation can
be precisely controlled by forcing.
\end{abstract}

\section{Introduction}

The automorphism tower of a group $G$ is obtained by iteratively
computing its automorphism group, the automorphism group of that
group, and so on transfinitely. Each group maps naturally into the
next via inner automorphisms and there is a natural direct limit
process (details below). 
$$G\to\Aut(G)\to\Aut(\Aut(G))\to\cdots\to G_\alpha\to G_{\alpha+1}\to\cdots$$
The tower {\it terminates} when a fixed point is first reached, a group that is isomorphic to its automorphism group by the natural map, and this
terminating ordinal is the {\it height} of the tower.

Although the automorphism tower construction has origins as a purely algebraic, group-theoretic construction, for some groups it has been observed to
exhibit an intriguing set theoretic behavior. For these groups, the automorphism tower is highly sensitive to the set-theoretic background in which
it is computed. For example, there can be a group whose automorphism tower is trivial in one model of set theory, but grows to uncountable heights in
other models of set theory, with all the intermediate heights also realized in still other models of set theory. Such set theoretic sensitivity is
both interesting and unusual for a purely algebraic construction.

\goodbreak

To give the details of the automorphism tower construction, one begins with any group $G_0=G$. At successor stages, we set
$G_{\alpha+1}=\Aut(G_\alpha)$, with the natural homomorphism $\pi_\alpha:G_\alpha\to G_{\alpha+1}$ sending a group element $g$ to the corresponding
inner automorphism $i_g:h\mapsto g^{-1}hg$, conjugating by $g$. At limit stages, $G_\lambda$ is the direct limit of the prior groups $G_\alpha$, for
$\alpha<\lambda$, with respect to the commutative system of homomorphisms $\pi_{\alpha,\beta}:G_\alpha\to G_\beta$ for $\alpha<\beta$ that are
obtained by composing the natural maps at each step. The tower terminates with height $\alpha$, if this is the earliest stage for which $\pi_\alpha$
is an isomorphism of $G_\alpha$ with $\Aut(G_\alpha)$. This occurs exactly when $G_\alpha$ is first a {\it complete} group, a centerless group having
only inner automorphisms. If the initial group $G$ is centerless, then all the groups in the tower are centerless, and so all the maps $\pi_\alpha$
are injective. In this special case, therefore, one may identify each group with its images in the later groups and thereby view the automorphism
tower as an increasing union of groups. Thomas \cite{Thomas:AutomorphismTowerProblem}, \cite{Thomas:AutomorphismTowerProblemII} proved that every
centerless group has a terminating automorphism tower. Building on this, Hamkins \cite{Hamkins98:EveryGroup} proved that every group leads eventually
to a centerless group, and consequently, every group has a terminating automorphism tower. Thomas' forthcoming monograph
\cite{Thomas:AutomorphismTowerProblemBook} is an excellent account of all aspects of the automorphism tower problem.

The possibility of groups whose automorphism towers are highly sensitive to set theory was established in \cite{HamkinsThomas2000:ChangingHeights},
where it was shown that there is a model of set theory, obtained by class forcing, in which the following statement holds:

\begin{statement}
\label{statement.HamkinsThomasChangingHeights}
For any ordinals $\alpha<\lambda$, there is a centerless group $G$ with an automorphism tower of height $\alpha$, but for every nonzero
$\beta<\lambda$, there is a cardinal and cofinality-preserving forcing extension in which the automorphism tower of the same group $G$ has height
$\beta$.
\end{statement}

To emphasize, \cite{HamkinsThomas2000:ChangingHeights} shows only that Statement \ref{statement.HamkinsThomasChangingHeights} is consistent, not that
it is a theorem; the groups fulfilling Statement \ref{statement.HamkinsThomasChangingHeights} are found in a forcing extension of the universe. It
was not known whether every model of \ZFC\ had such groups, or whether there were such groups, for example, in the constructible universe $L$. We
partially addressed this situation in our previous article \cite{FuchsHamkins:DegreesOfRigidity}, where we showed that $\diamondsuit$ implies that
there are groups whose automorphism towers are malleable by forcing, at least for the finite heights.
\begin{thm}\label{Theorem.DiamondImpliesGroups} (\cite{FuchsHamkins:DegreesOfRigidity})
Assume $\diamondsuit$ holds. Then for every $n<\omega$ there is a group $G$ with an automorphism tower of height $n$, but for any nonzero $m<\omega$
there is a forcing extension (preserving cardinals and cofinalities and not adding countable sequences of ordinals), in which the automorphism tower
of $G$ has height $m$.
\end{thm}

In this article, we prove fully that Statement
\ref{statement.HamkinsThomasChangingHeights} holds in $L$. This result
fulfills the suggestion at the conclusion of
\cite{FuchsHamkins:DegreesOfRigidity} that one might use a suitable
$\diamondsuit_{\kappa^+}$ hypothesis to carry out the construction
with $\kappa^+$-Souslin trees. 

\goodbreak

\begin{mainthm} 
\ \label{MainTheorem}
Statement \ref{statement.HamkinsThomasChangingHeights} holds in the constructible universe $L$.\break More generally, if there are unboundedly many
regular cardinals $\kappa$ for which $2^{<\kappa}=\kappa$ and
$\diamondsuit_{\kappa^+}(\CF_\kappa)$ holds, then Statement
\ref{statement.HamkinsThomasChangingHeights} holds.
\end{mainthm}
For any cardinal $\kappa$, we use $\CF_\kappa$ here to denote the set $\{\alpha<\kappa^+\st\cf(\alpha)=\kappa\}$, and
$\diamondsuit_{\kappa^+}(\CF_\kappa)$ is the assertion that there is a sequence $\vec D=\<D_\alpha\st\alpha\in\CF_\kappa>$ such that for any
$A\subset\kappa^+$ the set $\{\alpha\in\CF_\kappa\st A\intersect\alpha=D_\alpha\}$ is stationary in $\kappa^+$. In $L$, the hypotheses that
$2^{<\kappa}=\kappa$ and $\diamondsuit_{\kappa^+}(\CF_\kappa)$ are known to hold of every regular cardinal $\kappa$. Note that
$\diamondsuit_{\kappa^+}(\CF_\kappa)$ implies that $\kappa$ is regular, for otherwise $\CF_\kappa$ is empty.

The main theorem is proven by showing that it is possible from the
assumption to construct sequences of Souslin trees with a particular
combinatorial property known to imply Statement
\ref{statement.HamkinsThomasChangingHeights} (see Theorem
\ref{thm.CombinatorialCharacterization}). While we shoot directly
for the application to the automorphism tower problem in section
\ref{section:TheConstruction}, we show how to obtain sequences of
trees satisfying stronger combinatorial properties, which are of
independent interest, in section \ref{section:RealizingRelations}. 

Namely, let $\vT=\kla{T_\gamma\st\gamma<\lambda}$ be a sequence 
of rigid $\delta$-Souslin trees. Let $E$
be an equivalence relation on
$\lambda$. We say that $\vT$
\emph{realizes} $E$ if for all
$\alpha,\beta<\lambda$, $\alpha$ and $\beta$ are $E$-equivalent
iff $T_\alpha$ and
$T_\beta$ are isomorphic.
We say that $\vT$ is \emph{able to realize $E$} if
there is a cofinality-preserving, $<\!\lambda$-distributive notion of
forcing $\P_E$, such that in every $\P_E$-generic extension, $\vT$ is
a sequence of rigid Souslin trees which realizes $E$. By dropping the
condition that $\P_E$ be $<\!\lambda$-distributive, one arrives at the
concept of  sequences which are \emph{weakly able} to realize
equivalence relations. 
We call $E$ \emph{bounded} if there is some $\beta<\lambda$ such
that for all $\alpha\in(\beta,\lambda)$, $[\alpha]_E=\{\alpha\}$. 

The main result of section \ref{section:RealizingRelations}, Theorem
\ref{thm:RealizingSimpleEquivalenceRelations}, is:

\begin{thmStar}
{\sloppy Assume
$2^{<\kappa}=\kappa\ +\ \diamondsuit_{\kappa^+}(\CF_\kappa)$. Then there is 
a sequence $\langle T_\gamma$ $\st$ $\gamma<\kappa \rangle$ of $\kappa^+$-Souslin 
trees which is able to realize every bounded equivalence relation $E$
on $\kappa$.}
\end{thmStar}

The
result concerning weak realizability, Theorem \ref{thm.LongWeaklyRealizableSequences}, is:

\begin{thmStar}
Assume $2^{<\kappa}=\kappa\
+\ \diamondsuit_{\kappa^+}(\CF_\kappa)$. Then there is a sequence
$\langle T_\alpha\st\alpha<\kappa^+\rangle$ of $\kappa^+$-Souslin trees which
is weakly able to realize every equivalence relation on $\kappa^+$.
\end{thmStar}

\section{The Construction}
\label{section:TheConstruction}

The main algebraic construction of \cite{HamkinsThomas2000:ChangingHeights} shows that Statement \ref{statement.HamkinsThomasChangingHeights} is a
consequence of the existence of certain independent families of rigid graphs, whose isomorphism relation can be precisely controlled by forcing. The
malleable groups are then found within the automorphism group of a graph consisting of many disjoint copies of these graphs, and the corresponding
automorphism tower is controlled by forcing to control the isomorphism relation on these unit graphs. For the details of this construction, we refer
the reader to \cite{HamkinsThomas2000:ChangingHeights}, as well as to the overview article \cite{Hamkins2001:HowTall?} and to Thomas' excellent
forthcoming monograph \cite{Thomas:AutomorphismTowerProblemBook}.

In order to obtain the initial independent family of graphs, the construction of \cite{HamkinsThomas2000:ChangingHeights} uses forcing to add generic
Souslin trees, which exhibit the desired independence property, when construed as graphs, in the forcing extension. This is the reason that
\cite{HamkinsThomas2000:ChangingHeights} finds the malleable groups only in a forcing extension. The main contribution of this article is to replace
this forcing argument with a combinatorial construction; we replace the generic Souslin trees with Souslin trees constructed from a suitable diamond
hypothesis, the hypothesis that there are arbitrarily large cardinals $\kappa$ with $\diamondsuit_{\kappa^+}(\CF_\kappa)$ and $2^{<\kappa}=\kappa$.
Since this diamond hypothesis holds in the constructible universe $L$, we obtain both the trees and the malleable groups in $L$. So, let us state
without much further explanation that the construction of Section 2 in \cite{HamkinsThomas2000:ChangingHeights} shows that our Main Theorem is a
consequence of the following Theorem \ref{thm.CombinatorialCharacterization}, on which we shall now concentrate.

Before proceeding, let us be clear about our terminology. By a {\it $\kappa^+$-normal $\alpha$-tree}, we mean a tree of height $\alpha<\kappa^+$ with
all levels having size at most $\kappa$, such that there is a unique root node, every node (except those on the top level, if any) has at least two
immediate successors and has successors at all higher levels up to $\alpha$, and such that every node on a limit level is determined by its
predecessors. The tree is {\it $2$-splitting} if every node (except those on the top level, if any) has exactly two immediate successors. We write
$T|\alpha$ for the restriction of $T$ to levels less than $\alpha$, and we write $T(\gamma)$ for the $\gamma^{\rm th}$ level of $T$. An $\alpha$-tree
$T$ is ${<}\kappa$-closed if on every limit level $\delta<\alpha$ of cofinality less than $\kappa$, all cofinal branches through $T|\delta$ are
extended to $T(\delta)$. We write $[T]$ for the set of cofinal branches through any tree.  A notion of forcing $\P$ is {\it $\kappa$-distributive} if
it adds no new $\kappa$-sequences of elements of the ground model.

\goodbreak

\begin{thm}
\label{thm.CombinatorialCharacterization} Assume $\diamondsuit_{\kappa^+}(\CF_\kappa)$ and $2^{<\kappa}=\kappa$. Then there is a sequence
$\kla{T^\mu\st\mu<\kappa}$ of pairwise non-isomorphic rigid $\kappa^+$-Souslin trees such that:
\begin{enumerate}
\item
    \label{item.MakingTwoTreesIsomorphic}
    For every $\mu<\kappa$, there is a forcing extension in which $T^0$ becomes
    isomorphic to $T^\mu$, but the trees otherwise remain pairwise non-isomorphic and all remain rigid. The forcing furthermore preserves all
    cardinals and cofinalities and is $\kappa$-distributive.
\item
    \label{item.MakingEverythingBelowAlphaIsomorphic}
    For every $\mu<\kappa$, there is a forcing extension in which the all trees
    $T^\beta$ for $\beta<\mu$ become isomorphic, but the trees otherwise remain pairwise non-isomorphic and all remain rigid.
    The forcing furthermore preserves all cardinals and cofinalities and is $\kappa$-distributive.
\end{enumerate}
\end{thm}

\keinr \emph{Remark:} We shall show in section \ref{section:RealizingRelations} how to
handle other, more arbitrary patterns.

\proof We will construct the sequence $\vT=\kla{T^\mu\st\mu<\kappa}$ of $\kappa^+$-Souslin trees by simultaneous recursion on their levels, along
with a sequence of controller trees $\<C^\mu\st 0<\mu<\kappa>$, which will also be $\kappa^+$-Souslin trees and which as forcing notions will help us
to force the desired isomorphism patterns for $\vT$. Specifically, each $T^\mu$ will be a rigid $\kappa^+$ Souslin tree, and these trees will be
pairwise non-isomorphic. Each controller tree $C^\mu$ will also be a $\kappa^+$ Souslin tree, and when used as a notion of forcing, it will force
$T^0\cong T^\mu$ in such a way so as to fulfill statement \ref{item.MakingTwoTreesIsomorphic} of Theorem \ref{thm.CombinatorialCharacterization}. In
particular, forcing with $C^\mu$ will not create any unwanted isomorphisms between other pairs of trees, will preserve the rigidity of all of the
trees, and (being a ${<}\kappa$-closed $\kappa^+$ Souslin tree) will preserve all cardinals and cofinalities and be $\kappa$-distributive. The full
support product of controller trees $C^{<\mu}=\prod_{0<\beta<\mu}C^\beta$ will fulfill statement \ref{item.MakingEverythingBelowAlphaIsomorphic}. In
particular, this product will force that all $T^\beta$ for $\beta<\mu$ become isomorphic, that the trees otherwise remain pairwise non-isomorphic and
that all the trees remain rigid. By ensuring that $C^{<\mu}$ is actually a ${<}\kappa$-closed $\kappa^+$ Souslin tree, we will also ensure that it
preserves all cardinals and cofinalities and is $\kappa$-distributive. We use full support in this product because it is important for the
application to the Main Theorem that the forcing be $\kappa$-distributive; this conforms with and generalizes the usage of finite products of
controller trees in \cite{FuchsHamkins:DegreesOfRigidity}, where it was necessary to force only finitely many trees at a time to become isomorphic.
%

In order to achieve all these properties, we will use a diamond sequence from our hypothesis to anticipate and then
seal or kill off the various kinds of unwanted objects, such as unwanted uncountable antichains in the trees or
unwanted potential isomorphisms between the trees. It will be useful to have the following version of the diamond
hypothesis, allowing us more easily to anticipate such objects. Generalizing an idea from \cite{TSP}, for
$\alpha<\kappa^+$, let's let
$$H_{\kappa^+}(\alpha)=H_{\kappa^+}\cap V_\alpha.$$

\begin{lem}
\label{lem:ConsequenceOfDiamond} Assume $\diamondsuit_{\kappa^+}(\CF_\kappa)$ and $2^{<\kappa}=\kappa$. Then there is a
sequence\break $\<E_\alpha\st\alpha\in\CF_\kappa>$, with $E_\alpha\sub H_{\kappa^+}(\alpha)$, such that for every
$A\sub H_{\kappa^+}$, the set
    \[ \{\alpha\in\CF_\kappa\st
    A\cap H_{\kappa^+}(\alpha)=E_\alpha\}\] is stationary in $\kappa^+$.
\end{lem}

\proof Since $\card(H_{\kappa^+})=\kappa^+$, we may identify members of $H_{\kappa^+}$ with ordinals less than
$\kappa^+$. Let $f:H_{\kappa^+}\To\kappa^+$ be such an identification. Clearly, the set of $\alpha<\kappa^+$ such that
$f\rest H_{\kappa^+}(\alpha)$ is a bijection between $H_{\kappa^+}(\alpha)$ and $\alpha$, is club in $\kappa^+$. The
rest of the argument is a standard coding and decoding procedure. \qed

We now begin the detailed construction of the trees. Fix a diamond sequence $\<E_\alpha\st\alpha\in \CF_\kappa>$ as in Lemma
\ref{lem:ConsequenceOfDiamond}. All the trees $T^\mu$ and $C^\mu$ will be subtrees of ${}^{<\kappa^+}2$, ordered by initial segment. We want to
arrange that if $b$ is a generic branch through $C^\mu$, then there is an isomorphism $\pi_b$ of $T^0$ with $T^\mu$. The isomorphism $\pi_b$
witnessing this will be the map obtained by symmetric difference with $b$. Specifically, for any binary sequence $s$ let $\pi_s$ be the map swapping
the bits on any binary sequence at the coordinates that are $1$ in $s$. More precisely, $\pi_s(t)=s+t\!\mod2$, where we interpret the sum to have the
same length as $t$ (by padding $s$ with $0$s when $t$ is longer than $s$). This operation corresponds exactly to taking the symmetric difference up
to $|t|$ of the sets of which $s$ and $t$ are the characteristic function. Equivalently, we may define $|\pi_s(t)|=|t|$ and $\pi_s(t)(i)=t(i)$ if and
only if $s(i)=0$ or $i>|s|$. This collection of automorphisms has many convenient properties. For example, the maps commute because
$\pi_s\pi_t=\pi_{s+t}=\pi_{t+s}=\pi_t\pi_s$; they are all self-inverse, having order two because $\pi_s\pi_s=\pi_{s+s}=\pi_{\vec 0}=\id$; and they
have the convenient composition property that $\pi_{\pi_s(t)}=\pi_s\pi_t$, provided $|s|\geq|t|$. For any sequence $\vs=\<s_0,\ldots,s_{n-1}>$ of
binary sequences, we write $\pi_\vs$ for the composition $\pi_{s_0}\circ\cdots\circ\pi_{s_{n-1}}$.

To begin the construction, suppose that we have constructed the trees $T^\mu|\alpha$ and $C^\mu|\alpha$ below level $\alpha$. We will now define the
$\alpha^{\rm th}$ levels $T^\mu(\alpha)$ and $C^\mu(\alpha)$. We inductively assume that our trees satisfy the following conditions.

\begin{itemize}
\item[$(\star)_\alpha$]
\begin{enumerate}
 \item
   \label{item:Normality}
   Each $T^\mu|\alpha$ and $C^{\mu}|\alpha$ is a $2$-splitting
   $\kappa^+$-normal $\alpha$-tree as a subtree of
   ${}^{<\alpha}2$.
\item
  \label{item:Closedness}
  Each $T^\mu|\alpha$ and $C^{\mu}|\alpha$ is
  ${<}\kappa$-closed, in the tree order.
 \item
   \label{item:NodesAreIsomorphisms}
   If $\gamma<\alpha$ and $s\in C^\mu(\gamma)$, then
   $\pi_s\restrict T^0|(\gamma+1)$ is an isomorphism of
   $T^0|(\gamma+1)$ with $T^\mu|(\gamma+1)$.
%
\end{enumerate}
\end{itemize}%
%
%

The trees all begin, of course, with the empty root node $\<>$. Because of our insistence that the trees be $2$-splitting, we have no choice at
successor levels $\alpha+1$ but to extend every node on the $\alpha^{\rm th}$ level with its two immediate successors in ${}^{\alpha+1}2$. By doing
so, if the prior trees satisfy $(\star)_{\alpha+1}$, then it is easy to check that the resulting trees will satisfy $(\star)_{\alpha+2}$, and so we
will maintain our inductive assumption. At limit stages of the construction, if the conditions $(\star)_\alpha$ hold at all levels $\alpha$ below a
limit ordinal $\lambda$, then we automatically attain $(\star)_\lambda$ for the limit trees, because the $(\star)_\lambda$ hypothesis makes
assertions only about features of the trees occurring below level $\lambda$.

What remains is to construct the limit levels of the trees. We assume that the trees $T^\mu|\lambda$ and $C^\mu|\lambda$ are defined up to a limit
ordinal level $\lambda$ in such a way that $(\star)_\lambda$ is satisfied, and we must construct the $\lambda^{\rm th}$ levels of the trees
$T^\mu(\lambda)$ and $C^\mu(\lambda)$ in such a way that $(\star)_{\lambda+1}$ is satisfied. Defining the $\lambda^{\rm th}$ level of the trees
amounts to specifying for each tree the set of cofinal branches up to $\lambda$ which are to be extended. Since the trees consist of binary
sequences, we identify a branch $b$ with the binary sequence $\cup b$ extending it.

The easy limit case occurs when $\cf(\lambda)<\kappa$. In this case, in order to satisfy condition \ref{item:Closedness} of $(\star)_{\lambda+1}$, we
must extend every cofinal branch through every tree, defining $T^\mu(\lambda)=[T^\mu|\lambda]$ and $C^\mu(\lambda)=[C^\mu|\lambda]$. If
$\theta=\cf(\lambda)$, then since a branch is determined by its values on a cofinal set of levels and the earlier levels of the trees all have size
at most $\kappa$, the number of such branches in each case is at most $\kappa^\theta\leq \kappa^{{<}\kappa}=\kappa$, by our hypothesis. The extended
trees therefore remain normal and satisfy conditions \ref{item:Normality} and \ref{item:Closedness} of $(\star)_{\lambda+1}$. Condition
\ref{item:NodesAreIsomorphisms} is satisfied because if $b$ is a branch through $C^\mu|\lambda$ and $c$ is a branch through $T^0$, then $\pi_b''c$ is
a branch through $T^\mu$ and vice versa. Since all branches are extended, the nodes of $C^\mu$ at level $\lambda$ give rise to isomorphisms of $T^0$
with $T^\mu$. So our extended trees satisfy $(\star)_{\lambda+1}$, as desired.

The nontrivial limit case, the heart of our construction, occurs when $\cf(\lambda)=\kappa$. In this case, following
the general strategy of \cite{FuchsHamkins:DegreesOfRigidity}, we make use of the $\diamondsuit_{\kappa^+}(\CF_\kappa)$
sequence to seal various unwanted objects associated with the trees, such as unwanted maximal antichains, unwanted
automorphisms of the trees and unwanted isomorphisms between the trees. Before explicitly using the $\diamondsuit$
sequence, however, we shall separate the construction somewhat from the proof that it works, explaining in the Sealing
Lemma the sorts of unwanted objects that we can seal by selectively extending branches through the trees. The
terminology will be explained in the proof. For convenience, we write $\vec T|\lambda$ for the sequence
$\<T^\gamma|\lambda\st\gamma<\kappa>$ and $\vec C|\lambda$ for $\<C^\mu|\lambda\st 0<\mu<\kappa>$. We write
$C^{<\mu}|\lambda$ for the product $\prod_{0<i<\mu}(C^i|\lambda)$ with full support.

\goodbreak

\begin{sealinglemma}
\label{sublem:Sealing} Assume that $\lambda$ has cofinality $\kappa$ and that $(\star)_\lambda$ holds.
\begin{enumerate}
 \item \label{item:Extendibility}
$\vec
T|\lambda$ and $\vec C|\lambda$ can be extended to
$\vec{T}|(\lambda+1)$ and $\vec{C}|(\lambda+1)$ in such a way that
$(\star)_{\lambda+1}$ holds.
 \item
\label{item:SealingMaximalAntichainsInObjectTrees} If $A$ is a maximal antichain in $T^\gamma$, then $\vT|\lambda$ and $\vec{C}|\lambda$ can be
extended in such a way that $(\star)_{\lambda+1}$ holds and $A$ is sealed in $T^\gamma|(\lambda+1)$, meaning that every element of
$T^\gamma(\lambda)$ lies above an element of $A$.
 \item
\label{item:SealingMaximalAntichains} If $A$ is a maximal antichain in $C^{<\mu}|\lambda$, then $\vT|\lambda$ and $\vec{C}|\lambda$ can be extended
in such a way that $(\star)_{\lambda+1}$ holds and $A$ is sealed in $C^{<\mu}|(\lambda+1)$.
 \item
\label{item:SealingAutomorphismsOfTnu} If $f$ is a nontrivial automorphism of $T^\gamma|\lambda$, then $\vec T|\lambda$
and $\vec C|\lambda$ can be extended in such a way that $(\star)_{\lambda+1}$ holds and $f$ is sealed, meaning that $f$
cannot be extended to an automorphism of $T^\gamma|(\lambda+1)$.
 \item
\label{item:SealingCPotentialIsomorphisms} If $f$ is a $C^{<\mu}|\lambda$-potential isomorphism of $T^\gamma|\lambda$ and $T^\delta|\lambda$, where
$0<\mu\le\delta<\kappa$ and $\gamma<\delta$, then $\vec{T}|\lambda$ and $\vec{C}|\lambda$ can be extended in such a way that $(\star)_{\lambda+1}$
holds and $f$ is sealed, meaning that $f$ cannot be extended to a $C^{<\mu}|(\lambda+1)$-potential isomorphism of $T^\gamma|(\lambda+1)$ and
$T^\delta|(\lambda+1)$.
 \item
\label{item:SealingCPotentialAutomorphisms} If $f$ is a $C^{<\mu}|\lambda$-potential automorphism of $T^\gamma|\lambda$ then $\vec{T}|\lambda$ and
$\vec{C}|\lambda$ can be extended in such a way that $(\star)_{\lambda+1}$ holds and $f$ is sealed, meaning that $f$ cannot be extended to a
$C^{<\mu}|(\lambda+1)$-potential automorphism of $T^\gamma|(\lambda+1)$.
\end{enumerate}
\end{sealinglemma}

\proof We begin by describing our basic method for adding a $\lambda^{\rm th}$ level to the trees so as to ensure
$(\star)_{\lambda+1}$ for the extended trees. We will use this same construction template in all the subsequent cases.
To use this method, we first specify the $\lambda^{\rm th}$ level of the controller trees $C^\mu(\lambda)$ in such a
way that these continue to be $\kappa^+$-normal $(\lambda+1)$-trees. This amounts to choosing a covering set of
branches $C^\mu(\lambda)\subset [C^\mu|\lambda]$ of size at most $\kappa$ for each $\mu<\kappa$. Second, we select an
ordinal $\mu_0<\kappa$, and for the tree $T^{\mu_0}$ we specify a {\it generating} set $\Gamma^{\mu_0}$ of at most
$\kappa$ many branches covering $T^{\mu_0}|\lambda$. These branches generate others, through all the various
$T^\mu|\lambda$, by the application of appropriate compositions of the isomorphisms arising from elements of the
various $C^\mu(\lambda)$ we have just specified, and we must include these generated branches in order to maintain our
inductive assumption that branches through the controller trees give rise to isomorphisms between the object trees.
Specifically, let us define that a sequence $\vs=\<s_0,\ldots,s_n>$ is a \emph{trail} from $\zeta_0$ to $\zeta_{n+1}$
if there is a sequence $\<\zeta_0,\ldots,\zeta_{n+1}>$ of ordinals, called the {\it checkpoints} of the trail, such
that every other $\zeta_i$ is equal to $0$ and every other $\zeta_i$ is nonzero, and such that $s_i\in
C^{\max(\zeta_i,\zeta_{i+1})}(\lambda)$, for all $i\le n$. By $(\star)_\lambda$, each $\pi_{s_i}$ is an isomorphism of
$T^{\zeta_i}|\lambda$ with $T^{\zeta_{i+1}}|\lambda$, and so the full composition $\pi_\vs$ is an isomorphism of
$T^{\zeta_0}|\lambda$ with $T^{\zeta_{n+1}}|\lambda$. Since we want that $\pi_{\vec s}$ should be an isomorphism of the
extended trees, we define the $\lambda^{\rm th}$ level of $T^\gamma$ to consist of the corresponding set of generated
branches:
\[ T^\gamma(\lambda)=\{\pi_\vs(b)\st b\in\Gamma^{\mu_0}\tx{and}\vs\tx{is a trail
  leading from} \mu_0\tx{to}\gamma\}. \]
This completes the description of our construction template.

Let us show that as long as we follow this pattern, the extended trees will satisfy $(\star)_{\lambda+1}$, and so our
induction hypothesis will be maintained. To verify the normality of the extended trees, it suffices that the
$\lambda^{\rm th}$ levels of the trees cover the prior tree and have size at most $\kappa$ (the other points of
normality are easy to check). For the controller trees, the construction pattern explicitly called for $C^\mu(\lambda)$
to cover $C^\mu|\lambda$ and have size at most $\kappa$. For the object trees, observe first that there are at most
$\kappa$ many generated branches in $T^\gamma(\lambda)$, since there are at most $\kappa$ many trails $\vec s$ and at
most $\kappa$ many elements of $\Gamma^{\mu_0}$. To see that these branches cover $T^\gamma|\lambda$, suppose $p\in
T^\gamma|\lambda$. Pick any branches $s_0\in C^{\mu_0}(\lambda)$ and $s_1\in C^\gamma(\lambda)$. Then $\vs=\<s_0,s_1>$
is a trail leading from $\mu_0$ to $\nu$. Thus, $q=\pi_\vs(p)=\pi_\vs^{-1}(p)$ is in $T^{\mu_0}$, and since
$\Gamma^{\mu_0}$ covers $T^{\mu_0}|\lambda$, there is a branch $b\in\Gamma^{\mu_0}$ extending $q$. It follows that
$\pi_\vs(b)$ is a branch in $T^\nu(\lambda)$ extending $p$. So the extended trees are all normal, and we have fulfilled
condition \ref{item:Normality} of $(\star)_{\lambda+1}$. Condition \ref{item:Closedness}, asserting that the extended
trees are ${<}\kappa$-closed, is immediate because we are in the case $\cf(\lambda)=\kappa$. Condition
\ref{item:NodesAreIsomorphisms} will be satisfied because of the way we defined the generated branches in
$T^\gamma(\lambda)$. Specifically, for any $\gamma<\kappa$, we have:
\begin{itemize}
\item If $b\in T^0(\lambda)$ and $c\in C^\gamma(\lambda)$, then $\pi_c(b)\in T^\gamma(\lambda)$.

This is clear by the definition of $T^0(\lambda)$, since there is a trail $\vs$ leading from $\mu_0$ to $0$ and a
branch $\bar b\in\Gamma^{\mu_0}$ such that $b=\pi_\vs(\bar b)$. Thus, $\vu=\vs\verl\<c>$ is a trail leading from
$\mu_0$ to $\gamma$ and consequently, $\pi_c(b)=\pi_\vu(\bar b)\in T^\gamma(\lambda)$.
 \item If $d\in T^\gamma(\lambda)$ and $c\in C^\gamma(\lambda)$, then
  $\pi_c^{-1}(d)=\pi_c(d)\in T^0(\lambda)$.

Again, by the definition of $T^\gamma(\lambda)$,
  there is a trail $\vs$ leading from $\mu_0$ to $\gamma$ and a branch
  $\bar d\in\Gamma^{\mu_0}$ such that $d=\pi_\vs(\bar d)$. So $\vu=\vs\verl\<c>$ is a trail leading from $\mu_0$ to $0$, and
  consequently,
  $\pi_c(d)=\pi_c(\pi_\vs(\bar d))=\pi_\vu(\bar d)\in T^0(\lambda)$.
\end{itemize}%
Our construction template therefore ensures $(\star)_{\lambda+1}$ for the extended trees. We shall now use this method
to prove each of the statements of the Sealing Lemma.

\prooff{\ref{item:Extendibility}} By the construction template, we need only find, for every $\mu<\kappa$ and some
$\mu_0<\kappa$, sets of branches $C^\mu(\lambda)\of[C^\mu|\lambda]$ and $\Gamma^{\mu_0}\of[T^{\mu_0}|\lambda]$ of size
at most $\kappa$, covering their respective trees. Choose $\mu_0$ arbitrarily. Since the trees $C^\mu|\lambda$ and
$T^{\mu_0}|\lambda$ each have at most $\kappa$ many nodes, and each node can easily be extended to a cofinal branch
(using the ${<}\kappa$-closure of the tree and the fact that $\cf(\lambda)=\kappa$), we can easily construct the
desired covering sets of branches, and therefore, by the argument above, we attain $(\star)_{\lambda+1}$ for the
extended trees, as desired.

But in preparation for the later cases, let us explain in somewhat more elaborate detail a method for choosing the
covering sets of branches $C^\mu(\lambda)$ and $\Gamma^{\mu_0}$. We shall use a pseudo forcing construction with the
following partial order, with ${<}\kappa$ support in each factor:
\[
\P=(T^{\mu_0}|\lambda)^\kappa\times\prod_{0<\nu<\kappa}(C^\nu|\lambda)^\kappa,
\]
The idea is that the filter will provide, in the first factor, the generating branches $\Gamma^{\mu_0}$, and in the
second factor, the branches $C^\nu(\lambda)$ for each nonzero $\nu<\kappa$. We view conditions in $\P$ as pairs
$\<v,\vw>$, where $v:\kappa\To T^{\mu_0}|\lambda$ and $\vw=\kla{w_\nu\st
  0<\nu<\kappa}$, such that for all $\nu\in(0,\kappa)$,
$w_\nu:\kappa\To C^\nu|\lambda$. Because there is ${<}\kappa$-support, we have $w_\nu(i)=\<>$ for all but
less-than-$\kappa$-many $\nu$ and $i$, and $v(i)=\<>$ for all but less-than-$\kappa$-many $i$. A sufficiently generic
filter $H$ in $\P$ determines the generating branches in the first factor as follows:
\begin{ea*}
b_i&=&\bigcup\{v(i)\st\exists\vw\quad\<v,\vw>\in H\},\tx{for} i<\kappa,\\
\Gamma^{\mu_0}&=&\{b_i\st i<\kappa\}
\end{ea*}%
The controller branches are determined from $H$ by the second factor:
\begin{ea*}
c^{\nu}_i&=&\bigcup\{w_\nu(i)\st\exists v\quad \<v,\vw>\in H\},\tx{for}
i<\kappa,\ 0<\nu<\kappa,\\
C^\nu(\lambda)&=&\{c^\nu_i\st i<\kappa\}, \tx{for} 0<\nu<\kappa.
\end{ea*}%

We ensure that these sets of branches have the desired properties by ensuring that the filter $H$ meets certain dense
sets. Since the trees $T^{\mu_0}|\lambda$ and $C^{\nu}|\lambda$ are all ${<}\kappa$-closed, it follows that $\P$ is
${<}\kappa$-closed as a notion of forcing. A simple diagonalization then shows that, given any list of at most $\kappa$
many dense subsets of $\P$, there is a filter $H$ meeting each of them. We ensure that the $\Gamma^{\mu_0}$ and
$C^\nu(\lambda)$ arising from $H$ as above cover their respective trees by ensuring that $H$ meets the following dense
sets:
\begin{ea*}
\label{eqn:CoveringSets}
D_p&=&\{\<v,\vw>\in\P\st\exists i<\kappa\quad v(i)\ge p\},\tx{for every} p\in
T^{\mu_0}|\lambda,\\
D^\nu_q&=&\{\<v,\vw>\in\P\st\exists i<\kappa\quad w_\nu(i)\ge q\},\tx{for
  every} \nu\in(0,\kappa)\\
&& \hskip2in\tx{and every} p\in T^{\mu_0}|\lambda.\nonumber
\end{ea*}%
We ensure that the branches $b_i$ and $c^\nu_i$ are cofinal by meeting the dense sets:
\begin{ea*}
D_{\alpha,i}&=&\{\<v,\vw>\st |v(i)|>\alpha\},\tx{for}
\alpha<\lambda \tx{and} i<\kappa,\\
D_{\nu,\alpha,i}&=&\{\<v,\vw>\st |w_\nu(i)|>\alpha\}, \tx{for}
0<\nu<\kappa,\ \alpha<\lambda \tx{and} i<\kappa.
\end{ea*}%
Altogether, we have $\kappa$ many dense sets, so there is a filter $H$ meeting them all. Thus, we have constructed covering sets of cofinal branches
$\Gamma^{\mu_0}$ and $C^\nu(\lambda)$, as desired, and so statement \ref{item:Extendibility} is proved.\qedd{(\ref{item:Extendibility})}

\prooff{\ref{item:SealingMaximalAntichainsInObjectTrees}} Suppose that $A$ is a maximal antichain in $T^\gamma|\lambda$. We shall follow the
construction template, taking $\mu_0=\gamma$. First, we specify the covering sets $C^\nu(\lambda)$ arbitrarily. We now carry out a pseudo forcing
construction with the poset $\Q=(T^\gamma|\lambda)^\kappa$, using ${<}\kappa$-support. A filter $H$ in $\Q$ will add the desired generating branches
$\Gamma^\gamma$, just as in the first factor of $\P$ above. Meeting the dense sets $D_p =\{v\in\Q\st\exists i<\kappa\quad v(i)\ge p\}$ for $p\in
T^\gamma|\lambda$ and $D_{\alpha,i}=\{ v\in\Q\st |v(i)|>\alpha\}$ for $\alpha<\lambda$ and $i<\kappa$ ensures that the resulting $\Gamma^\gamma$
covers $T^\gamma|\lambda$ with cofinal branches. It remains to ensure that $A$ is sealed. For any $i<\kappa$ and any trail $\vt$ from $\gamma$ to
$\gamma$, consider the following dense set.
\[ D_{\vt,i}=\{v\in\Q\st\exists p\in A\quad\pi_\vt(v(i))\ge p\} \]
To see that it is dense, suppose $v\in\Q$. Let $p=v(i)$ and $q=\pi_\vt(p)$. By the maximality of $A$, we can extend $q$ to some $q'$ above an element
of $A$. It follows that $p'=\pi_\vt^{-1}(q')$, which is the same as $\pi_\vt(q')$, is above $p$, because $\pi_\vt$ is an automorphism of
$T^\gamma|\lambda$. So, if $v'$ extends $v$ by extending the $i^{\rm th}$ coordinate from $p$ to $p'$, it follows that $v'\in D_{\vt,i}$, and so it
is dense. Finally, if $H$ meets all the $D_{\vt,i}$, then $A$ will be sealed in $T^\gamma|(\gamma+1)$, because these dense sets exactly ensure that
every element of $T^\gamma(\lambda)$ will lay above an element of $A$. So we have sealed $A$ while retaining $(\star)_{\lambda+1}$.
\qedd{(\ref{item:SealingMaximalAntichainsInObjectTrees})}

\prooff{\ref{item:SealingMaximalAntichains}} Next, we seal maximal antichains in the controller product trees
$C^{<\mu}|\lambda$. Suppose that $A$ is a maximal antichain in $C^{<\mu}|\lambda$. We shall build our covering sets of
branches $C^{\nu}(\lambda)$ so that nodes $\vec b=\<b_\nu\st 0<\nu<\mu>\in\prod_{0<\nu<\mu}C^\nu(\lambda)$ all lie
above a node in $A$. This is what we mean by sealing the antichain. We use the partial order $\P$ as above, with
$\mu_0$ chosen arbitrarily, and construct a pseudo generic filter $H$ by meeting a list of $\kappa$ many dense sets. We
can ensure that the sets of branches $\Gamma^{\mu_0}$ and $C^\nu(\lambda)$ resulting from $H$ cover their respective
trees and consist of cofinal branches by meeting the dense sets mentioned in the proof of statement
\ref{item:Extendibility}. To ensure that $A$ is sealed, we construct $H$ to meet the following dense sets, for every
$\vec l\in\kappa^\mu$:
\begin{ea*}
D_{\vec l}&=&\{\<v,\vw>\in\P\st\exists r\in A\quad \<w_\nu(l_\nu)\st\nu<\mu>\ge r\}
\end{ea*}%
Each of these sets is dense, precisely because $A$ is a maximal antichain. So we can extend the trees in such a way that $A$ is sealed and
$(\star)_{\lambda+1}$ holds.\qedd{(\ref{item:SealingMaximalAntichains})}

\prooff{\ref{item:SealingAutomorphismsOfTnu}} Although statement \ref{item:SealingAutomorphismsOfTnu} is a consequence
of statement \ref{item:SealingCPotentialAutomorphisms}, we prove this easier case first in order to introduce the
technique in a less complicated situation. Suppose that $f$ is a nontrivial automorphism of the tree
$T^{\mu_0}|\lambda$, taking $\mu_0$ to be $\gamma$ of the statement of \ref{item:SealingAutomorphismsOfTnu}. We want to
extend the trees in such a way that $f$ does not extend to an automorphism of $T^{\mu_0}|(\lambda+1)$, while retaining
$(\star)_{\lambda+1}$.
%
%
To do so, we shall construct a pseudo generic filter $H$ in the partial order $\P$ as above. In addition to meeting the dense sets of statement
\ref{item:Extendibility}, which ensure that the branches provided by $H$ cover the trees and are cofinal, we shall ensure that $f$ is sealed by
meeting additional dense sets.

Specifically, in order to seal $f$ we will arrange that there is a fixed generating branch $b_0\in\Gamma^{\mu_0}$ such that whenever $\vs$ is a trail
leading from $\mu_0$ to $\mu_0$, and $b$ is a generating branch in $\Gamma^{\mu_0}$, then $\pi_{\vs}(b)\neq f[b_0]$. This will seal $f$, because
$b_0$ is in $T^{\mu_0}(\lambda)$, but according to the construction template, $f[b_0]$ will not be added to $T^{\mu_0}(\lambda)$, and so $f$ will not
extend to an automorphism of $T^{\mu_0}|(\lambda+1)$. Note that we allow $\vs$ to be the empty trail, interpreting $\pi_\vs$ in this case as the
identity function. The difficulty, of course, is that we don't know the trails leading from $\mu_0$ to $\mu_0$ before specifying the sets
$C^\nu(\lambda)$. But we do know how such trails will arise via $\P$. Let us therefore define that a \emph{template} for a trail, from $\zeta_0$ to
$\zeta_{n+1}$, is a pair $\ttt=\kla{\<i_0,\ldots,i_n>,\<\zeta_0,\ldots,\zeta_{n+1}>}$ such that every other $\zeta_l$ is zero, every other $\zeta_l$
is non-zero, and each $\zeta_l$ and $i_l$ is less than $\kappa$. The idea is that $\ttt$ is a template for the trail
\[\<c^{\max(\zeta_0,\zeta_1)}_{i_0},\ldots,c^{\max(\zeta_n,\zeta_{n+1})}_{i_n}>\]
that will ultimately be determined by the filter $H$. Any condition $\<v,\vw>\in\P$ gives partial information about these branches $c^\zeta_i$ and
consequently also partial information about this trail, which we denote:
\[ \ttt_\kla{v,\vw}=\<w_{\max(\zeta_0,\zeta_1)}(i_0),\ldots,w_{\max(\zeta_n,\zeta_{n+1})}(i_n)>.\]

We now describe the dense sets that will ensure that $f$ is sealed. Since $f$ is a nontrivial automorphism of $T^{\mu_0}|\lambda$, there is a node
$p_0\in T^{\mu_0}|\lambda$ that is moved by $f$. Suppose that $\ttt$ is a template for a trail leading from $\mu_0$ to $\mu_0$, and let $i<\kappa$.
Let $\hat p$ be the condition in $\P$ placing $p_0$ onto the first generating branch $b_0$. That is, $\hat{p}=\kla{v,\vw}$ where $v(0)=p_0$ and
otherwise $v(j)=\<>$ and $w_\nu(j)=\<>$. We claim that the following set is dense in $\P$ below $\hat p$.
\[ D_{f,\ttt,i}=\{u=\kla{v,\vw}\in\P\st
\pi_{\ttt_u}(v(i))\perp f(v(0))\}.
\]
To see that this set is dense below $\hat p$, suppose $u=\<v,\vw>$ is any condition in $\P$ below $\hat p$. If $\ttt$ is the empty template, or more
generally if the template has internal cancellation causing $\pi_{\ttt_u}$ to necessarily be the identity function, then it is relatively easy to
extend $u$ to a condition in $D_{f,\ttt,i}$, using the fact that $f(p_0)\perp p_0$ in the case $i=0$. So suppose that
$\ttt=\kla{\<i_0,\ldots,i_n>,\<\zeta_0,\ldots,\zeta_{n+1}>}$ corresponds to a non-trivial $\pi_\ttt$. Because the maps $\pi_s$ all have order two and
commute, it follows that at least one of the branches $c^{\max(\zeta_k,\zeta_{k+1})}_{i_k}$ specified by the trail $\ttt$ appears an odd number of
times in $\ttt$. The first step is to extend $u$ to a condition $u'$ that specifies the partial information about the branches through the controller
trees relevant for the computation of $\ttt$ and also the coordinate $0$ of $v$ all to the same height. The next step is to extend $u'$ to $u'_0$ and
$u'_1$, which extend these branches one bit further, in an identical way, except that $u'_0$ and $u'_1$ differ on this extra bit for that odd branch
coordinate $c^{\max(\zeta_k,\zeta_{k+1})}_{i_k}$. It follows that $\pi_{\ttt_{u'_0}}(v(0))\neq \pi_{\ttt_{u'_1}}(v(0))$, and so one of them must be
incompatible with $f(v(0))$. Thus, either $u'_0$ or $u'_1$ is in $D_{f,\ttt,i}$ below $u$, and the set is dense. Since there are at most $\kappa$
many templates for trails, we may construct a pseudo generic filter $H$ in $\P$ meeting all the dense sets we have mentioned. It follows that the
resulting sets of branches $\Gamma^{\mu_0}$ and $C^\nu(\lambda)$ cover their respective trees, consist of cofinal branches and seal $f$, as
desired.\qedd{(\ref{item:SealingAutomorphismsOfTnu})}

\prooff{\ref{item:SealingCPotentialIsomorphisms}} Following terminology introduced in
\cite{FuchsHamkins:DegreesOfRigidity}, we say that $f$ is a {\it $C^{<\mu}|\lambda$-potential isomorphism} of
$T^\gamma|\lambda$ with $T^\delta|\lambda$, if $f$ is an order preserving function from $C^{<\mu}|\lambda$ into the
collection of partial isomorphisms of $T^{\gamma}|\lambda$ to $T^{\delta}|\lambda$, such that for any condition $\vec
q\in C^{<\mu}|\lambda$, there is a dense set of stronger conditions $\vec r$, whose $f(\vec r)$ extends $f(\vec q)$ so
as to insert any given node of $T^{\gamma}|\lambda$ into the domain and any given node of $T^{\delta}|\lambda$ into the
range. Such functions $f$ arise exactly from $C^{<\mu}$-names for isomorphisms of $T^{\gamma}$ to $T^{\delta}$, as one
may take $f(\vec q)$ as the information about that name forced by $\vec q$. This is merely a relatively concrete way to
treat such names.

So suppose that $f$ is a $C^{<\mu}|\lambda$-potential isomorphism of $T^\gamma|\lambda$ with $T^\delta|\lambda$, where
$0<\mu\le\delta<\kappa$ and $\gamma<\delta$. We shall extend the trees to level $\lambda$ by constructing a pseudo
generic filter $H$ in the partial order $\P$, using $\mu_0=\gamma$, and applying the construction template. As before,
we may ensure that the sets of branches $\Gamma^{\mu_0}$ and $C^\nu(\lambda)$ arising from $H$ cover their respective
trees and consist of cofinal branches, by meeting at most $\kappa$ many dense sets. In order to ensure also that $f$ is
sealed, we now specify some additional dense sets.

To explain our strategy for sealing $f$, let us imagine for a moment that $H$ has been already selected, giving rise to
the covering sets $\Gamma^{\mu_0}$ and $C^\nu(\lambda)$, for $0<\nu<\kappa$. We would like to have branches $c_\nu\in
C^\nu(\lambda)$, for $0<\nu<\mu$, and a generating branch $b\in\Gamma^{\mu_0}$, such that $f[\vc]$ is an isomorphism
from $T^{\gamma}|\lambda$ to $T^{\delta}|\lambda$, but such that for every trail $\vt$ leading from $\gamma$ to
$\delta$ and every generating branch $d\in\Gamma^{\mu_0}$ we have $f[\vc][b]\neq\pi_\vt(d)$. This expresses precisely
that $f[\vc][b]$ is not one of the generated branches constituting $T^{\delta}(\lambda)$. If we can accomplish this,
then $f$ will not extend to a potential isomorphism of the extended trees, since the partial isomorphism $f[\vc]$ will
not extend to an isomorphism that works on level $\lambda$, and so $f$ will be sealed. We will set things up in such a
way that if $H$ is generic with respect to the dense sets we specify, then the witnessing branches $c_\nu$ for the
above strategy will be the branches $c^\nu_0$, as defined from $H$, and the branch $b$ will be $b_0$, as defined from
$H$, using the notation for the branches as in the proof of statement \ref{item:Extendibility} above.

First, in order to ensure that $f[\vc]$ is an isomorphism between $T^{\gamma}|\lambda$ and $T^{\delta}|\lambda$, it suffices that $H$ intersect each
of the following subsets of $\P$. These sets are dense precisely because $f$ is a potential isomorphism of $T^\gamma|\lambda$ with
$T^\delta|\lambda$.
\begin{ea*}
D^0_{f,p}&=&\{\kla{v,\vw}\st p\in\dom(f(\<w_i(0)\st i<\mu>))\},
\tx{for every} p\in T^\gamma,\\
D^1_{f,q}&=&\{\kla{v,\vw}\st q\in\ran(f(\<w_i(0)\st i<\mu>))\}, \tx{for every} q\in T^\delta.
\end{ea*}%
Next, for each template $\ttt$ for a trail leading from $\gamma$ to $\delta$ and each $i<\kappa$, we will have $H$ intersect the following dense set.
\begin{ea*}
D_{f,\ttt,i}=\{u=\<v,\vw>\in\P\st f(\<w_\xi(0)\st 0<\xi<\mu>)(v(0))\perp\pi_{\ttt_u}(v(i))\}
\end{ea*}%
This set ensures that the strategy we mentioned above will be realized, because conditions in it exactly ensure that
$f[\vc][b]\neq\pi_\vt(d)$, using $\vc=\<c^\xi_0\st 0<\xi<\mu>$ and the generating branches $b=b_0$ and $d=b_i$, as we
explained above. To see that $D_{f,\ttt,i}$ is dense, we make critical use of the fact that $\mu\leq\delta$. Given any
condition $u\in \P$, we first extend its $w_\xi$'s for $\xi<\mu$ and its $v(0)$ and $v(i)$ so that $v(0)$ is in the
domain of the part of $f$ ``decided'' by it, and so that $v(i)$ is at the same height as $v(0)$, which is larger than
the height of the coordinates of the controller trees specified by $u$ that occur in the template trail $\ttt$. Now
there must be such controller tree coordinates, since the template trail leads to $\delta$ (so $C^\delta$ is involved),
which is at least $\mu$, while in order to ``decide'' $f$, only controller coordinates below $\mu$ are needed (since
$f$ has domain $C^{{<}\mu}|\lambda$). Fixing all but one such coordinate, and then extending the remaining one in
different ways will result in conditions $u'$ with different outcomes for $\pi_{\ttt_{u'}}(v(i))$. One of these
outcomes must therefore be different from $f(\<w_\xi(0)\st 0<\xi<\mu>)(v(0))$, and so the resulting condition $u'$ will
be in $D_{f,\ttt,i}$, showing that it is dense.

In summary, if $H$ meets all the dense sets we have mentioned, then we will have successfully accomplished our strategy for extending the trees in
such a way that $f$ is sealed and $(\star)_{\lambda+1}$ holds.\qedd{(\ref{item:SealingCPotentialIsomorphisms})}

\goodbreak

\prooff{\ref{item:SealingCPotentialAutomorphisms}} Suppose that $f$ is a $C^{<\mu}|\lambda$-potential automorphism of $T^\gamma|\lambda$. Following a
strategy similar to that in case \ref{item:SealingCPotentialIsomorphisms}, we will again specify a collection of dense subsets of $\P$, using
$\mu_0=\gamma$, such that any pseudo generic filter $H$ meeting them will give rise to the desired tree extensions according to the construction
template. As above, with $\kappa$ many dense sets we can easily ensure that the branch sets $\Gamma^{\mu_0}$ and $C^\nu(\lambda)$ arising from $H$ do
indeed cover their respective trees and consist of cofinal branches.

To explain our strategy for sealing $f$, let us again imagine that $H$ has already been chosen. We will arrange that there is a sequence
$\vc=\<c_\nu\st 0<\nu<\mu>$ of controller branches with $c_\nu\in C^\nu(\lambda)$ and a branch $b\in\Gamma^{\gamma}$, such that $f[\vc]$ is an
automorphism of $T^\gamma|\lambda$, but such that for any trail $\vt$ leading from $\gamma$ to $\gamma$ and every generating branch
$d\in\Gamma^\gamma$ we have $f[\vc](b)\neq\pi_\vt(d)$. This strategy will seal $f$, because we will have added $\vc$ to the controller product
$C^{{<}\mu}|(\lambda+1)$, but $f[\vc]$ will not extend to an automorphism of $T^\gamma|(\lambda+1)$, because $b$ is a branch there, while $f[\vc][b]$
is not. To carry out this strategy, it will suffice that $H$ meet certain dense sets, which force that the controller branches $c_\nu=c^\nu_0$ and
generating branch $b=b_0$ will witness the desired property.

First, in order to ensure that $f[\vb]$ is an automorphism of $T^\gamma|\lambda$, it suffices that $H$ intersects the following dense sets:
\begin{ea*}
D^0_{f,p}&=&\{\kla{v,\vw}\st p\in\dom(f(\<w_i(0)\st i<\mu>))\},
\tx{for every} p\in T^\gamma,\\
D^1_{f,q}&=&\{\kla{v,\vw}\st q\in\ran(f(\<w_i(0)\st i<\mu>))\},
\tx{for every} q\in T^\gamma.
\end{ea*}%
We may work below a condition $\<v,\vw>$ such that $f(\<w_\nu(0)\st 0<\nu<\mu>)(v(0))\perp v(0)$, which will allow us to realize $b$ (above) as
$b_0$.

Next, to fulfill the second part of our strategy, suppose that $\ttt$ is a template for a trail leading from $\gamma$ to $\gamma$ and that
$i<\kappa$, we will have $H$ intersect the following set.
\[
D_{f,\ttt,i}=\{u=\<v,\vw>\in\P\st f(\<w_\xi(0)\st
0<\xi<\mu>)(v(0))\perp\pi_{\ttt_u}(v(i))\}.
\]%
This set ensures that our strategy will be realized, because conditions in it exactly ensure that
$f[\vc](b)\neq\pi_\vt(d)$, using $\vc=\<c^\xi_0\st 0<\xi<\mu>$ and the generating branches $b=b_0$ and $d=b_i$, as we
explained above. It remains only to check that $D_{f,\ttt,i}$ is dense. This is clear when $0<i$, since $v(i)$ can be
extended in incompatible ways, giving rise to different values of $\pi_{\ttt_u}(v(i))$ with the same value of $f(\vec
w)(v(0))$, causing one of the extensions to be in $D_{f,\ttt,i}$. So we may assume that $i=0$. Suppose
$\ttt=\kla{\<i_0,\ldots,i_n>,\<\zeta_0,\ldots,\zeta_{n+1}>}$. As in case \ref{item:SealingCPotentialIsomorphisms}, if
$\ttt$ is trivial is the sense that it gives rise only to the identity function $\pi_{\ttt_u}$, then it is easy to
extend a condition into $D_{f,\ttt,i}$ using the fact that we are working under the condition $\<v,\vw>$ forcing that
$f(\vw)(v(0)\perp v(0)$. So we may assume that $\ttt$ is nontrivial. It follows, using the fact that the maps $\pi_s$
all commute and have order two, that one of the coordinate pairs $\<\max(\zeta_k,\zeta_{k+1}),i_k>$ appearing in $\ttt$
appears an odd number of times. Because $\ttt$ is a template for a trail from $\gamma$ to $\gamma$, a ``closed'' trail
if you will, it has an even number of coordinate pairs altogether, and more precisely, every checkpoint $\nu$ used in
$\ttt$ is used twice each time it appears, once going from $\nu$ to $0$ and once from $0$ to $\nu$. It follows that the
trail cannot always use branch index $0$ in each of these directions, that is, not every coordinate pair in $\ttt$ has
the form $\<\nu,0>$, corresponding to the branch $c^\nu_0$, because in this case the branches $c^\nu_0$ would all be
used an even number of times, causing them to cancel and make $\ttt$ trivial after all. Thus, there must be some
coordinate pair in $\ttt$ of the form $\<\nu,i>$ with $i\neq 0$. This allows us to argue as in case
\ref{item:SealingCPotentialIsomorphisms}, by specifying everything but this one coordinate sufficiently high, and then
considering two incompatible extensions of this one coordinate. More precisely, given any condition $\<v',\vw'>$, we
extend only its $w'_\nu(i)$ in incompatible ways. This leads to incompatible outcomes on the right hand side of the
formula defining $D_{f,\ttt,i}$, while the left hand side is the same. So one of these possibilities of extending
$\<v',\vw'>$  must yield a different outcome on the right hand side than on the left, and we have found a stronger
condition in $D_{f,\ttt,i}$, thereby verifying that this set is
dense.\qedd{(\ref{item:SealingCPotentialAutomorphisms})}

The proof of the Sealing Lemma is now complete.\qed 

We now continue with the proof of Theorem \ref{thm.CombinatorialCharacterization} and the recursive construction of the
sequences $\vec{T}$ and $\vec{C}$. To remind the reader of our context, we are in the case where $\lambda$ is a limit
ordinal of cofinality $\kappa$, and we have constructed the trees $\vec{T}|\lambda$ and $\vec{C}|\lambda$ in such a way
that $(\star)_\lambda$ holds. We have also fixed a $\diamondsuit_{\kappa^+}(\CF_\kappa)$ sequence $\vec{E}$
anticipating subsets of ${}^{<\kappa}(H_{\kappa^+})$ in the sense of Lemma \ref{lem:ConsequenceOfDiamond}, of which we
shall now make critical use. We now extend the trees to level $\lambda$ in such a way to attain $(\star)_{\lambda+1}$
for the extended trees, while also working to seal various unwanted objects, by dividing into cases depending on the
value of the diamond sequence $E_\lambda$.

\emph{Case 1.} If $E_\lambda=\{\<1,\gamma,s>\st s\in A\}$, where $A$ is a maximal antichain in $T^{\gamma}|\lambda$,
then we extend the trees so as to seal the antichain, according to the Sealing Lemma statement
\ref{item:SealingMaximalAntichainsInObjectTrees}.

\emph{Case 2.} If $E_\lambda=\{\<2,\mu,\vx>\st\vx\in A\}$, for some maximal antichain $A$ in $C^{<\mu}|\lambda$, then
we extend the trees so as to seal the antichain, according to the Sealing Lemma statement
\ref{item:SealingMaximalAntichains}.

\emph{Case 3.} If $E_\lambda=\{\<3,\gamma,s,t>\st f(s)=t\}$, for some $f$ which is an automorphism of
$T^\gamma|\lambda$, then we extend $\vec{T}|\lambda$ and $\vec{C}|\lambda$ according to the Sealing Lemma statement
\ref{item:SealingAutomorphismsOfTnu}, thereby sealing $f$ as an automorphism of $T^\gamma|\lambda$.

\emph{Case 4.} If $E_\lambda=\{\<4,\mu,\gamma,\delta,s,t,\vx>\st\vx\in C^{<\mu}|\lambda\tx{and} f(\vx)(s)=t\}$, where
$f$ is a $C^{<\mu}|\lambda$-potential isomorphism between $T^{\gamma}|\lambda$ and $T^\delta|\lambda$ and
$0<\mu\le\delta<\kappa$ and $\gamma<\delta$, then we extend $\vec{T}|\lambda$ and $\vec{C}|\lambda$ according to the
Sealing Lemma statement \ref{item:SealingCPotentialIsomorphisms}, so as to seal $f$ as a $C^{<\mu}|\lambda$-potential
isomorphism automorphism between $T^\gamma|\lambda$ and $T^\delta|\lambda$.

\emph{Case 5.} If $E_\lambda=\{\<5,\mu,\gamma,s,t,\vx>\st\vx\in C^{<\mu}|\lambda\tx{and} f(\vx)(s)=t\}$, where $f$ is a
$C^{<\mu}|\lambda$-potential automorphism of $T^\gamma|\lambda$, then we extend $\vec{T}|\lambda$ and $\vec{C}|\lambda$
according to the Sealing Lemma statement \ref{item:SealingCPotentialAutomorphisms}, thereby sealing $f$ as a
$C^{<\mu}|\lambda$-potential automorphism of between $T^\gamma|\lambda$.

\emph{Case 6.} Finally, if none of the above cases occur, then we use the Sealing Lemma statement
\ref{item:Extendibility} to extend the trees in such a way that $(\star)_{\lambda+1}$ holds for the extended trees.

\goodbreak

This completes the recursive construction of the trees $T^\gamma$ and the controller trees $C^\mu$. We now prove that
the trees are as we claimed. We observe first that all the trees are in fact Souslin trees. It is clear that $T^\gamma$
and $C^\mu$ are $\kappa^+$-normal trees of height $\kappa^+$, because the recursive construction proceeds in such a way
that the restrictions $T^\gamma|\lambda$ and $C^\mu|\lambda$ are $\kappa^+$-normal $\lambda$-trees for every
$\lambda<\kappa^+$. The object tree $T^\gamma$ is now seen to be Souslin by the usual reflective argument. Namely, if
$A\of T^\gamma$ is any maximal antichain, then there is a closed unbounded set of stages $\alpha$ such that
$A\intersect T^\gamma|\alpha$ is a maximal antichain in $T^\gamma|\alpha$. Since the set $A'=\{\<1,\mu,s>\st s\in A\}$
is anticipated by the diamond sequence $\vE$ on a stationary set, there will be stationarily many stages
$\lambda\in\CF_\kappa$, such that $E_\lambda=A'\intersect H_{\kappa^+}(\lambda)$ and $A\intersect T^\gamma|\lambda$ is
a maximal antichain in $T^\gamma|\lambda$, putting us in case 1 of the construction. At such a stage, the tree
$T^\gamma|(\lambda+1)$ was specifically designed to seal $A$. All elements in the tree above level $\lambda$,
consequently, are compatible with an element of $A\intersect T^\gamma|\lambda$, and so $A\subset T^\gamma|\lambda$.
Thus, $A$ has size at most $\kappa$, and so $T^\gamma$ is in fact a Souslin tree. An essentially identical argument
shows that the controller product tree $C^{{<}\mu}$ is also a $\kappa^+$-normal $\kappa^+$-Souslin tree, using case 2.
Specifically, if $A\of C^{{<}\mu}$ is a maximal antichain, then the coding set $A'=\{\<2,\mu,\vx>\st\vx\in A\}$ is
anticipated by $\vE$ on a stationary subset of $\CF_\kappa$, and so there is a stage $\lambda\in\CF_\kappa$ such that
$E_\lambda=A'\intersect H_{\kappa^+}(\lambda)$, for which $A\intersect C^{{<}\mu}|\lambda$ is a maximal antichain. By
the construction in case 2, therefore, the antichain $A$ is sealed, and consequently is contained in
$C^{{<}\mu}|\lambda$, which has size at most $\kappa$. So $C^{{<}\mu}$ is a Souslin tree, as desired. It follows
directly from this that the individual controller trees $C^\nu$ are also $\kappa^+$-Souslin trees.

Second, we observe that the controller trees create the desired isomorphisms. By design, any branch $s$ through the
controller tree $C^{\mu}$ creates an isomorphism $\pi_s$ from $T^0$ to $T^\mu$. Forcing with $C^\mu$, therefore, will
ensure $T^0\cong T^\mu$. Since $C^\mu$ is a ${<}\kappa$-closed $\kappa^+$-Souslin tree, it preserves all cardinals and
cofinalities and is $\kappa$-distributive. Similarly, forcing with $C^{{<}\mu}$ will force all the object trees $T^\nu$
for $\nu<\mu$ to be isomorphic, because the product forcing adds generic branches through every individual factor
$C^\nu$ for $0<\nu<\mu$. And since the controller product $C^{{<}\mu}$ is a ${<}\kappa$-closed $\kappa^+$-Souslin tree,
it preserves all cardinals and cofinalities, and is $\kappa$-distributive.

Next, we check that the controller trees preserve the rigidity of the object trees. Suppose towards contradiction that
forcing with the controller product $C^{<\mu}$ created a non-trivial automorphism of some object tree $T^\gamma$. Then
there would be a condition $\vec p\in C^{<\mu}$ and a name $\dot{\pi}$ such that $\vec p$ forces via $C^{<\mu}$ that
$\dot{\pi}$ is a nontrivial automorphism of $T^\gamma$. Let $f$ be the function mapping any $\vec q$ extending $\vec p$
in $C^{<\mu}$ to the part of $\dot{\pi}$ that is decided by $\vec q$. That is, $f(\vec q)=\{\<s,t>\st \vec
q\forces\dot{\pi}(\check s)=\check t\}$. If $\vec q$ is incompatible with $\vec p$, then let $f(\vec q)$ be the
identity function. Using the fact that $C^{{<}\mu}$ is $\kappa$-distributive, as we established above, it follows that
all the proper initial segments of the automorphism named by $\dot\pi$ are in the ground model, and so $f$ is a
$C^{{<}\mu}$-potential automorphism of $T^\gamma$. An easy argument shows that there is therefore a club set of stages
$\lambda$ for which $f\restrict(C^{<\mu}|\lambda)$ is a $C^{{<}\mu}|\lambda$-potential automorphism of
$T^\gamma|\lambda$. By coding this potential automorphism and using the diamond sequence to anticipate it via case 4,
it follows that at some such stage $\lambda$ we sealed this potential automorphism as in the Sealing Lemma statement
\ref{item:SealingCPotentialAutomorphisms}, a contradiction since $f\restrict(C^{{<}\mu}|\lambda)$ extends to $f$. So
forcing with $C^{{<}\mu}$ preserves the rigidity of all the object trees $T^\gamma$, as deired. It now follows directly
that forcing with just one controller tree $C^\nu$ will also preserve the rigidity of all these trees, since $C^\nu$
appears as a factor in $C^{{<}\mu}$ for any $\mu>\nu$. It also follows that the trees $T^\mu$ will all be rigid, since
any actual non-trivial automorphism can easily be used to construct a potential automorphism.

Finally, we observe that an essentially similar argument shows that the controller trees and controller products do not
create unwanted isomorphisms between the object trees, using case 4 of the construction and Sealing Lemma statement
\ref{item:SealingCPotentialIsomorphisms}. From this, it also follows that the object trees $T^\gamma$ are pairwise
non-isomorphic.

Thus, the proof of Theorem \ref{thm.CombinatorialCharacterization} is now complete. \qed

\section{Realizing equivalence relations}
\label{section:RealizingRelations}

Let us briefly review what we have done in section \ref{section:TheConstruction} of this paper. We proved
Theorem \ref{thm.CombinatorialCharacterization}, under the assumption of $2^{<\kappa}+\diamondsuit_{\kappa^+}(\CF_\kappa)$, and
deduced from this, using the main algebraic construction from \cite{HamkinsThomas2000:ChangingHeights}, that Statement
\ref{statement.HamkinsThomasChangingHeights} is true for $\lambda=\kappa^+$.

But actually, Theorem \ref{thm.CombinatorialCharacterization} is weaker and less natural than the combinatorial criterion given
in \cite{HamkinsThomas2000:ChangingHeights}. It suffices for our
application, as can be verified by looking at the
algebraic construction, but since the original statement is of independent interest, we restate it here and ask whether
we can construct Souslin trees having properties along the lines of
the original combinatorial criterion. 
The combinatorial property shown consistent by forcing in \cite{HamkinsThomas2000:ChangingHeights} is that for any
regular cardinal $\lambda$, there is a sequence of $\lambda$-Souslin
trees which is able to realize every equivalence relation on $\lambda$, see the
introduction for the relevant definitions.

What we get by refining the construction from section \ref{section:TheConstruction} is:

\begin{thm}
\label{thm:RealizingSimpleEquivalenceRelations} Assume $2^{<\kappa}=\kappa+\diamondsuit_{\kappa^+}(\CF_\kappa)$. Then there is
a sequence $\kla{T_\gamma\st\gamma<\kappa}$ of $\kappa^+$-Souslin
trees which is able to realize every bounded equivalence relation on
$\kappa$.
\end{thm}

Again, the notion of a bounded equivalence relation was defined in the
introduction.
If we want full realizability, we seem forced to climb up in
cardinality one further step:

\begin{cor}
\label{cor:RealizingEquivalenceRelations} Assume $2^\kappa=\kappa^++\diamondsuit_{\kappa^{++}}(\CF_{\kappa^+})$. Then
there is a sequence $\kla{T_\alpha\st\alpha<\kappa}$ of
$\kappa^{++}$-Souslin trees that is able to realize every equivalence
relation on $\kappa$.
\end{cor}

\prooff{Theorem \ref{thm:RealizingSimpleEquivalenceRelations}}
(Sketch). 

The construction produces controller trees
$C^{\mu,\nu}$, for $\mu<\nu<\kappa$, which are intended to add an isomorphism between $T^\mu$ and $T^\nu$. Towards
realizing such an equivalence relation $E$, let
\[I_E=\{\kla{\mu,\nu}\st\mu<\nu\ \tx{and} \mu\tx{is least such that} \mu
E\nu\}.\]
The aim is that the forcing realizing $E$ is
\[ C_E=\prod_{\kla{\mu,\nu}\in I_E}C^{\mu,\nu}.\]
This is the reason for the restriction to bounded equivalence relations:
We want this product to consist of less than $\kappa$ many components;
otherwise we wouldn't be able to anticipate names for objects we want to seal via the
$\diamondsuit_{\kappa^+}(\CF_\kappa)$ sequence.

Notice the similarity between every single component $\prod_{\nu E\mu, \nu>\mu}C^{\mu,\nu}$ (with some fixed
$\mu=\min[\mu]_E$) here and the product $C^{<\mu}$ that we worked with before. The additional complication in the
present situation is that there may be many such components in the product. Let's look a little more closely at the
details of the construction. The heart of the construction is again the (modified version of the) Sealing Lemma
\ref{sublem:Sealing}. Assuming we have constructed the $\vT|\lambda$ and $\vC|\lambda$, the critical points are:
\begin{enumerate}
\item
\label{item:SealingCPotentialIsomorphismsAgain} If $f$ is a $C_E|\lambda$-potential isomorphism of $T^\gamma|\lambda$
and $T^\delta|\lambda$, where $\gamma\not\!\!{E} \delta$, then $\vec{T}|\lambda$ and $\vec{C}|\lambda$ can be extended
in such a way that $(\star)_{\lambda+1}$ holds and $f$ is sealed.
\item
\label{item:SealingCPotentialAutomorphismsAgain}
If $f$ is a $C_E|\lambda$-potential automorphism of $T^\gamma|\lambda$
then $\vec{T}|\lambda$ and $\vec{C}|\lambda$ can be extended in such a way that $(\star)_{\lambda+1}$ holds and $f$ is
sealed.
\end{enumerate}
Here, $(\star)_{\lambda+1}$ is the obvious analog of what we worked with before.

Here is a sketch of the proof of \ref{item:SealingCPotentialIsomorphismsAgain}, in the context of Theorem
\ref{thm:RealizingSimpleEquivalenceRelations}. The construction template is as before, \emph{mutatis mutandis}. Thus,
we specify the $\lambda^{\rm th}$ level of the controller trees $C^{\mu,\nu}(\lambda)$ in such a way that these
continue to be $\kappa^+$-normal $(\lambda+1)$-trees. Also, we select an ordinal $\mu_0<\kappa$, and for the tree
$T^{\mu_0}$, we specify a {\it generating} set $\Gamma^{\mu_0}$ of at most $\kappa$ many branches covering
$T^{\mu_0}|\lambda$. The branches of the trees in $\vT|\lambda$ that are going to be extended will then be those
generated by $\Gamma^{\mu_0}$ under trail embeddings. This time, a sequence $\vs=\kla{s_0,\ldots,s_n}$ is a
\emph{trail} from $\zeta_0$ to $\zeta_{n+1}$ if there is a sequence $\<\zeta_0,\ldots,\zeta_{n+1}>$ of ordinals, called
the {\it checkpoints} of the trail, such that $s_i\in C^{\zeta_i,\zeta_{i+1}}(\lambda)$, for all $i\le n$. Here, we use
the notation $C^{\mu,\nu}=C^{\nu,\mu}$.

In analogy to the previous construction, the method for determining the covers of the controller trees and the
generating set of branches through $T^{\mu_0}$ is by a pseudo forcing construction with the following partial order,
with ${<}\kappa$ support in each factor:
\[
\P=(T^{\mu_0}|\lambda)^\kappa\times\prod_{\mu<\nu<\kappa}(C^{\mu,\nu}|\lambda)^\kappa.
\]
We view conditions in $\P$ as pairs $\<v,\vw>$, where $v:\kappa\To T^{\mu_0}|\lambda$ and
$\vw=\kla{w_{\mu,\nu}\st\mu<\nu<\kappa}$, such that $w_{\mu,\nu}:\kappa\To C^{\mu,\nu}|\lambda$. If $H$ is sufficiently
$\P$-generic, we set:
\begin{ea*}
b_i&=&\bigcup\{v(i)\st\exists\vw\quad\<v,\vw>\in H\},\tx{for} i<\kappa,\\
\Gamma^{\mu_0}&=&\{b_i\st i<\kappa\},\\
c^{\mu,\nu}_i&=&\bigcup\{w_{\mu,\nu}(i)\st\exists v\quad \<v,\vw>\in H\},\tx{for}
i<\kappa,\ 0<\nu<\kappa,\\
C^{\mu,\nu}(\lambda)&=&\{c^{\mu,\nu}_i\st i<\kappa\}, \tx{for} \mu<\nu<\kappa.
\end{ea*}%

So suppose that $f$ is a $C_E|\lambda$-potential isomorphism of $T^\gamma|\lambda$ with $T^\delta|\lambda$, where
$\gamma$ and $\delta$ are not $E$-equivalent. Set $\mu_0=\gamma$. As before, we may ensure that the sets of branches
$\Gamma^{\mu_0}$ and $C^\nu(\lambda)$ arising from $H$ cover their respective trees and consist of cofinal branches, by
meeting certain dense sets in $\P$. In order to ensure also that $f$ is sealed, we now specify some additional dense
sets.

We want to choose $H$ in such a way that the corresponding covering sets $\Gamma^{\mu_0}$ and $C^{\mu,\nu}(\lambda)$,
for $\mu<\nu<\kappa$ have the property that there are branches $c_{\mu,\nu}\in C^{\mu,\nu}(\lambda)$, for
$\kla{\mu,\nu}$ in $I_E$, and a generating branch $b\in\Gamma^{\mu_0}$, such that $f[\vc]$ is an isomorphism from
$T^{\gamma}|\lambda$ to $T^{\delta}|\lambda$, but such that for every trail $\vt$ leading from $\gamma$ to $\delta$ and
every generating branch $d\in\Gamma^{\mu_0}$ we have $f[\vc][b]\neq\pi_\vt(d)$. This way, $f$ will not extend to a
potential isomorphism of the extended trees, since the partial isomorphism $f[\vc]$ will not extend to an isomorphism
that works on level $\lambda$, and so $f$ will be sealed. We will set things up in such a way that if $H$ is generic
with respect to the dense sets we specify, then the witnessing branches $c_{\mu,\nu}$ for the above strategy will be
the branches $c^{\mu,\nu}_0$, as defined from $H$, and the branch $b$ will be $b_0$, as defined from $H$, using the
notation for the branches as above.

Ensuring that $f[\vc]$ is an isomorphism between $T^{\gamma}|\lambda$ and $T^{\delta}|\lambda$ works as before.

Now, the modified notion of a trail comes with an analogous modification of the notion of a \emph{template for a
trail}. For each such template $\ttt$ for a trail leading from $\gamma$ to $\delta$ and each $i<\kappa$, we will have
$H$ intersect the following dense set:
\begin{ea*}
D_{f,\ttt,i}=\{u=\<v,\vw>\in\P\st f(\<w_{\xi,\zeta}(0)\st\kla{\xi,\zeta}\in I_E>)(v(0))\perp\pi_{\ttt_u}(v(i))\}.
\end{ea*}%
To see that $D_{f,\ttt,i}$ is dense, we make critical use of the fact that $\delta$ is not $E$-equivalent to $\mu$.
Given any condition $u\in \P$, we first extend its $w_{\xi,\zeta}$'s for $\kla{\xi,\zeta}\in I_E$, and its $v(0)$ and
$v(i)$ so that $v(0)$ is in the domain of the part of $f$ ``decided'' by it, and so that $v(i)$ is at the same height
as $v(0)$, which is larger than the height of the coordinates of the controller trees specified by $u$ that occur in
the template trail $\ttt$. Obviously, the trail template cannot be trivial, since it leads from $\gamma$ to $\delta$.
Moreover, there must be a coordinate $\kla{\kla{\mu,\nu},j}$ occuring in the trail template, such that $\mu$ and $\nu$
are not $E$-equivalent, and such that the coordinate doesn't cancel. Now $\kla{\mu,\nu}\notin I_E$, so this coordinate
is not needed in order to ``decide'' $f$. Fixing all but one such coordinate, and then extending it in different ways
will result in conditions $u'$ with different outcomes for $\pi_{\ttt_{u'}}(v(i))$. One of these outcomes must
therefore be different from $f(\<w_{\kla{\xi,\zeta}}(0)\st \kla{\xi,\zeta}\in I_E>)(v(0))$, and so the resulting
condition $u'$ will be in $D_{f,\ttt,i}$, showing that it is dense.

Now let's give a sketch of the proof of \ref{item:SealingCPotentialAutomorphismsAgain}.

Suppose that $f$ is a $C_E|\lambda$-potential automorphism of $T^\gamma|\lambda$. Following a strategy similar to that
in case \ref{item:SealingCPotentialIsomorphismsAgain}, we will again specify a collection of dense subsets of $\P$,
using $\mu_0=\gamma$, such that any pseudo generic filter $H$ meeting them will give rise to the desired tree
extensions according to the construction template. As above, with $\kappa$ many dense sets we can easily ensure that
the branch sets $\Gamma^{\mu_0}$ and $C^\nu(\lambda)$ arising from $H$ do indeed cover their respective trees and
consist of cofinal branches.

For the moment, let us again imagine that $H$ has already been chosen. We will arrange that there is a sequence
$\vc=\<c_{\xi,\zeta}\st \kla{\xi,\zeta}\in I_E>$ of controller branches with $c_{\xi,\zeta}\in C^{\xi,\zeta}(\lambda)$
and a branch $b\in\Gamma^{\gamma}$, such that $f[\vc]$ is an automorphism of $T^\gamma|\lambda$, but such that for any
trail $\vt$ leading from $\gamma$ to $\gamma$ and every generating branch $d\in\Gamma^\gamma$ we have
$f[\vc](b)\neq\pi_\vt(d)$. This strategy will seal $f$, because we will have added $\vc$ to the controller product
$C^{{<}\mu}|(\lambda+1)$, but $f[\vc]$ will not extend to an automorphism of $T^\gamma|(\lambda+1)$, because $b$ is a
branch there, while $f[\vc][b]$ is not. To carry out this strategy, it will suffice that $H$ meet certain dense sets,
which force that the controller branches $c_{\mu,\nu}=c^{\mu,\nu}_0$ and generating branch $b=b_0$ will witness the
desired property.

Ensuring that $f[\vc]$ is an automorphism of $T^\gamma|\lambda$ works as before. We may work below a condition
$\<v,\vw>$ such that $f(\<w_{\xi,\zeta}(0)\st \kla{\xi,\zeta}\in I_E>)(v(0))\perp v(0)$, which will allow us to realize
$b$ (above) as $b_0$.

Next, suppose that $\ttt$ is a template for a trail leading from $\gamma$ to $\gamma$ and that $i<\kappa$. We will have
$H$ intersect the following set.
\[
D_{f,\ttt,i}=\{u=\<v,\vw>\in\P\st f(\<w_{\xi,\zeta}(0)\st\kla{\xi,\zeta}\in I_E>)(v(0))\perp\pi_{\ttt_u}(v(i))\}.
\]%
It remains to check that $D_{f,\ttt,i}$ is dense. This is clear when $0<i$, since $v(i)$ can be extended in
incompatible ways, giving rise to different values of $\pi_{\ttt_u}(v(i))$ with the same value of $f(\vec w)(v(0))$,
causing one of the extensions to be in $D_{f,\ttt,i}$. So we may assume that $i=0$. Suppose
$\ttt=\kla{\<\zeta_0,\zeta_1,i_0>,\ldots,\<\zeta_n,\zeta_{n+1},i_n>}$. As in case
\ref{item:SealingCPotentialIsomorphismsAgain}, if $\ttt$ is trivial in the sense that it gives rise only to the
identity function $\pi_{\ttt_u}$, then it is easy to extend a condition into $D_{f,\ttt,i}$ using the fact that we are
working below the condition $\<v,\vw>$ forcing that $f((\vw\rest I_E)(0))(v(0))\perp v(0)$. So we may assume that
$\ttt$ is nontrivial. It follows, using the fact that the maps $\pi_s$ all commute and have order two, that one of the
triples $\<\zeta_k,\zeta_{k+1},i_k>$ appearing in $\ttt$ appears an odd number of times. Here, of course we have to
identify $\<\xi,\zeta,j>$ and $\<\zeta,\xi,j>$ when counting.

In order to run the argument that worked before, we have to find a coordinate in the trail template which is irrelevant
for the value of $f$. This is a little more involved in the current situation.

First, if there is a coordinate in the trail which does not stay within $I_E$, meaning that there is a coordinate
$\kla{\zeta_k,\zeta_{k+1},i_k}$ occurring in $\ttt$ such that neither $\kla{\zeta_k,\zeta_{k+1}}$ nor
$\kla{\zeta_{k+1},\zeta_k}$ is in $I_E$, and if this coordinate
doesn't cancel (i.e., if it occurs an odd number of times in
the trail, under the above identification), then we have found a coordinate with the desired properties, since $f$ only
depends on coordinates in $I_E$.

So now assume that every trail coordinate leaving $I_E$ cancels. Then the trail can be viewed as consisting of a series
of closed trails each of which stays within $I_E$ (if there are several closed component trails which stay in the same
equivalence class, then they may be viewed as one. Since the
isomorphisms we are dealing with commute, the order in which a trail
is hiked is irrelevant). One of these closed component trails is
nontrivial, or else the entire trail would have been trivial. Let's fix
some such trail. We can now apply the original argument. Again, remember the
similarity between every single component $\prod_{\nu E\mu,
  \nu>\mu}C^{\mu,\nu}$, with fixed
$\mu=\min[\mu]_E$ here and the product $C^{<\mu}$ that we worked with
before. The checkpoints of the trail we fixed bounce back and forth between
the minimum of the equivalence class within which it is staying and other
members of that equivalence class. So it cannot be the case that every
coordinate triple in that trail 
appears with branch index $0$, or else it would end up being trivial,
since the trail is closed. So there must be some coordinate triple that
doesn't cancel and has nonzero branch index. This is a coordinate that doesn't affect the value of $f$ but that is
relevant for the value of $\pi_{\ttt_u}(v(0))$.

Now we can argue as in case
\ref{item:SealingCPotentialAutomorphismsAgain}, by specifying
everything but this one
coordinate sufficiently high, and then considering two incompatible
extensions of this one coordinate. This leads to 
incompatible outcomes on the right hand side of the formula defining
$D_{f,\ttt,i}$, while the left hand side is the
the same. So one of these possibilities of extending must yield a
different outcome on the right hand side than on the
left, and we have found a stronger condition in $D_{f,\ttt,i}$,
thereby verifying that this set is dense.

This finishes the proof sketch of
\ref{item:SealingCPotentialAutomorphismsAgain}, and thus of the analog
of the Sealing Lemma \ref{sublem:Sealing} in the context of realizing
equivalence relations.

The construction now works as before. Basically, fixing a
$\diamond_{\kappa^+}(\CF_\kappa)$-sequence, it tells us what we have
to seal at stage $\lambda$ of the construction ($\lambda$ being a limit
ordinal of cofinality $\kappa$), and the (new version of the) Sealing
Lemma tells us that we can do that. The new objects we have to worry
about are $C_E|\lambda$-potential isomorphisms of $E$-inequivalent
trees $T_\gamma$ and $T_\delta$, and $C_E|\lambda$-potential
automorphisms. These instances of sealing correspond to cases 4 and 5
of the construction in Theorem
\ref{thm.CombinatorialCharacterization}. It is a simple matter to code
the equivalence relations in some canonical way into the elements of
the diamond sequence. Remember that they basically are subsets of
$\kappa$ (they can even be coded as bounded ones). The same thing has
to be done when sealing maximal antichains in $C_E|\lambda$, as in
case 2. Sealing antichains on the $\vT$-side works as before.

The argument that this construction achieves what we wanted, works as
before. Given some equivalence relation $E$ on $\kappa$, $C_E$ is a
Souslin tree, because if it had a maximal antichain of size
$\kappa^+$, this would reflect down to a stage of the construction
where it was sealed. So $C_E$ is $<\!\kappa^+$-distributive and
cofinality-preserving. In generic extensions by $C_E$, there are
generic branches in $C^{\alpha,\beta}$ whenever $\alpha$ and $\beta$
are $E$-equivalent and $\alpha$ is least with $\alpha E \beta$. These
give rise to isomorphisms between the 
trees $T^\alpha$ and $T^\beta$. By composing these isomorphisms, one
sees that $T^\gamma$ and $T^\delta$ are isomorphic in the extension
whenever $\gamma E\delta$.
If $\alpha$ and $\beta$ are
$E$-inequivalent, then no such isomorphism is added, or else there
would be a $C_E$-potential additional isomorphism between $T^\alpha$
and $T^\beta$. Again, this isomorphism would reflect down to a stage
of the construction where it was sealed in case 4. Similarly, the trees
$T^\alpha$ remain rigid because a $C_E$-potential automorphism of some
$T^\alpha$ would reflect down to a stage where it was sealed in case 5
of the construction.

This finishes the proof of Theorem \ref{thm:RealizingSimpleEquivalenceRelations}.\qed

\medskip

Now let's turn to weak realizability, which omits the requirement of $<\!\lambda$-distributivity from full
realizability, as defined in the introduction. We get:

\begin{thm}
\label{thm.LongWeaklyRealizableSequences} Assume $2^{<\kappa}=\kappa\ +\ \diamondsuit_{\kappa^+}(\CF_\kappa)$. Then
there is a sequence $\kla{T_\alpha\st\alpha<\kappa^+}$ of
$\kappa^+$-Souslin trees which is almost able to realize every
equivalence relation on its members.
\end{thm}

\proof (Sketch) We carry out the construction of Theorem \ref{thm:RealizingSimpleEquivalenceRelations}, but with longer
sequences of trees. We shall construct the Souslin trees $\<T^\gamma\st\gamma<\kappa^+>$ along with the controller
trees $\<C^{\mu,\nu}\st \mu<\nu<\kappa^+>$ by simultaneous recursion
on their levels. The aim is that in the end,  if we are given an
equivalence relation $E$ on $\kappa^+$, define $I_E$ as before and let
$C_E$ be the corresponding product of controller trees, but this time
with $<\!\kappa$-support, $C_E$ will witness that $\vT$ is weakly able to
realize $E$.

In order to achieve this, the following observation is useful.

\begin{lem}
Under the assumption of $2^{<\kappa}=\kappa\ +\ \diamondsuit_{\kappa^+}(\CF_\kappa)$, the sequence
$\kla{E_\alpha\st\alpha\in\CF_\kappa}$ of Lemma
\ref{lem:ConsequenceOfDiamond} in fact has the following property:

Whenever $A\sub{}^{<\kappa}H_{\kappa^+}$, the set
\[ \{\alpha\in\CF_\kappa\st
    A\cap{}^{<\kappa}(H_{\kappa^+}(\alpha))=E_\alpha\}\]
is stationary in $\kappa^+$.
\end{lem} 
\proof In contrast to the proof of \ref{lem:ConsequenceOfDiamond}, it is important here that the
$\diamondsuit_{\kappa^+}$-sequence is based on $\CF_\kappa$. Namely, for $\alpha<\kappa^+$ of cofinality $\kappa$, it
is the case that ${}^{<\kappa}H_{\kappa^+}(\alpha)\sub H_{\kappa^+}(\alpha)$. This is all that's needed to carry out
the proof. \qed

This allows us to anticipate both equivalence relations $E$ and most
of the other unwanted objects associated to forcing with $C_E$.

By anticipating
antichains, one ensures that $C_E$ satisfies the $\kappa^+$-c.c.; and since the controller trees are all
$<\!\kappa$-closed, so is $C_E$.

In the part of the argument where we ensure that forcing with $C_E$ doesn't add unwanted isomorphisms, however, we
are faced with a new problem. Namely, we aren't able to properly anticipate the isomorphisms that might be added. Since
the controller product forcing is not $\kappa$-distributive, it could happen that for some unwanted isomorphism added
by that forcing, the restriction of it to a level of the tree is not in the ground model. Thus, the diamond sequence
would not be able to anticipate and seal it. How can we handle such
potential isomorphisms? Our solution is the following 
trick: we will ensure a stronger rigidity property that can be anticipated.

More precisely, suppose we want to ensure that forcing with $C_E$ doesn't add an isomorphism between $T^\gamma$
and $T^\delta$, where $\gamma$ and $\delta$ are not $E$-equivalent. To
do this, we will instead ensure that forcing with $C_E\times 
T^{\gamma}$ doesn't add a branch to $T^\delta$. This amounts to anticipating and sealing a $C_E\times
T^\gamma$-potential additional branch of $T^\delta$, to use the terminology of \cite{FuchsHamkins:DegreesOfRigidity}.

At the same time, we have to anticipate the equivalence relation
itself. 

The construction relies on the following observation: If $E$ is an equivalence relation on $\kappa^+$,
and $f$ is, say, a $C_E\times T_\gamma$-potential additional
branch of $T_\delta$, where $\gamma$ and $\delta$ are not $E$-equivalent, then the set $C$ consisting of all
$\alpha<\kappa^+$ such that $f\rest((C_{E\rest\alpha}|\alpha)\times T_\gamma)$ is a $(C_{E\rest\alpha}|\alpha)\times
T_\gamma$-potential additional branch of $T_\delta$ is club in $\kappa^+$. By sealing such objects during the
construction, i.e., by anticipating both the potential additional branches and the equivalence relation, it is ensured
that there won't be such potential additional branches of the whole trees.

Similarly, we can ensure that $C_E$ preserves the rigidity of the trees by ensuring that forcing with
$C_E\times T^\gamma$ adds exactly one branch to $T^\gamma$. Thus, the trees $T^\gamma$ will have a strong
version of the unique branch property, in the terminology of \cite{FuchsHamkins:DegreesOfRigidity}. For this, we
anticipate and seal the $C^{{<}\mu}\cross T^\gamma$-potential additional branches of $T^\gamma$.

This finishes the sketch of the proof of Theorem \ref{thm.LongWeaklyRealizableSequences}.\qed

The previous construction raises the following question:

\begin{question}
In the constructible universe $L$, is there a cardinal $\lambda$ and a
$\lambda$-sequence of $\lambda$-Souslin trees that is able to
realize every equivalence relation on $\lambda$?
\end{question}

\bibliographystyle{hep.bst}
\bibliography{HamkinsBiblio,MathBiblio,literatur}

\begin{thebibliography}{Ham01}

\bibitem[DJ74]{TSP}
K.~J. Devlin and H.~Johnsbr{\aa}ten,
\newblock \textsl{ The Souslin Problem},
\newblock Lecture Notes in Mathematics {\bf 405}, Springer, Berlin, 1974.

\bibitem[FH06]{FuchsHamkins:DegreesOfRigidity}
G.~Fuchs and J.~D. Hamkins, \textsl{ {Degrees of rigidity for {S}ouslin
  trees}},
\newblock ArXiv Mathematics e-prints  (February 2006), {math.LO/0602482},
\newblock Submitted to the Journal of Symbolic Logic.

\bibitem[Ham98]{Hamkins98:EveryGroup}
J.~D. Hamkins, \textsl{ Every group has a terminating transfinite automorphism
  tower},
\newblock Proc. Amer. Math. Soc. \textbf{ 126}(11), 3223--3226 (1998),
  {arXiv:math.GR/9808014}.

\bibitem[Ham01]{Hamkins2001:HowTall?}
J.~D. Hamkins, \textsl{ How tall is the automorphism tower of a group?},
\newblock Logic and Algebra, AMS Contemporary Mathematics Series \textbf{ 302},
  49--57 (2001).

\bibitem[HT00]{HamkinsThomas2000:ChangingHeights}
J.~D. Hamkins and S.~Thomas, \textsl{ Changing the heights of automorphism
  towers},
\newblock Ann. Pure Appl. Logic \textbf{ 102}(1-2), 139--157 (2000),
  {arXiv:math.LO/9703204}.

\bibitem[Tho]{Thomas:AutomorphismTowerProblemBook}
S.~Thomas,
\newblock \textsl{ The Automorphism Tower Problem},
\newblock to appear.

\bibitem[Tho85]{Thomas:AutomorphismTowerProblem}
S.~Thomas, \textsl{ The automorphism tower problem},
\newblock Proceedings of the American Mathematical Society \textbf{ 95},
  166--168 (1985).

\bibitem[Tho98]{Thomas:AutomorphismTowerProblemII}
S.~Thomas, \textsl{ The automorphism tower problem {II}},
\newblock Israel Journal of Mathematics \textbf{ 103}, 93--109 (1998).

\end{thebibliography}

\end{document}